\documentclass[10pt,letterpaper]{article}
\usepackage{blindtext,titlefoot}
\usepackage{authblk}
\usepackage{marvosym}
\usepackage[utf8]{inputenc}
\usepackage[english]{babel}
\usepackage{amsmath}
\usepackage{amsfonts}
\usepackage{amssymb}
\usepackage{makeidx}
\usepackage{graphicx}
\usepackage{lmodern}
\usepackage{kpfonts}
\usepackage{dsfont}
\usepackage{cancel}
\usepackage{amsthm} 
\usepackage[left=2cm,right=2cm,top=2cm,bottom=2cm]{geometry}
\usepackage{subcaption}
\usepackage{multirow}
\usepackage{float} 
\usepackage{url}
\usepackage{breakurl}
\usepackage[breaklinks]{hyperref}
\usepackage{multicol}
\usepackage{mathtools, cuted}
\usepackage{lipsum}
\usepackage{booktabs}
\usepackage{comment}
\usepackage{cite}

\usepackage{longtable}
\DeclareMathOperator{\Op}{O}
\DeclareMathOperator{\Ma}{M}
\DeclareMathOperator{\IMa}{InvM}
\DeclareMathOperator{\Conv}{Conv}
\DeclareMathOperator{\Ord}{Ord}
\DeclareMathOperator{\ord}{ord}
\providecommand{\nset}[1]{
\mathbb{#1}
}
\providecommand{\set}[1]{
\left\{#1\right\}
}

\providecommand{\ifr}[5]{
{}^{#1}_{#2}{#3}_{#4}^{#5}
}

\providecommand{\gam}[1]{
\Gamma\left(#1 \right)
}

\providecommand{\norm}[1]{
\left\lVert #1 \right\rVert
}
\providecommand{\abs}[1]{
\left\lvert #1 \right\rvert
}
\providecommand{\ds}[1]{
\displaystyle #1
}
\providecommand{\der}[3]{
\dfrac{#1^{#3} }{ #1 #2^{#3}}
}

\newcommand{\mref}[2]{\textbf{#1 #2}}
\newtheorem{theorem}{ Theorem}[section]

\newtheorem{corollary}[theorem]{Corollary}
\newtheorem{example}{Example}

\usepackage{enumitem} 
\setlist[itemize]{noitemsep} 

\usepackage{titlesec} 
\titleformat{\section}[block]{\large\bfseries\scshape\centering}{\thesection.}{1em}{} 
\titleformat{\subsection}[block]{\large\bfseries\scshape\centering}{\thesubsection.}{1em}{}
\titleformat{\subsubsection}[block]{\large\bfseries\scshape\centering}{\thesubsubsection.}{1em}{} 

\title{\huge \bfseries Acceleration of the order of convergence of a family of fractional fixed point methods and its implementation in the solution of a nonlinear algebraic system related to hybrid solar receivers}

\author[,a]{A. Torres-Hernandez  \footnote{Email: anthony.torres@ciencias.unam.mx; Corresponding author; ORCID: 0000-0001-6496-9505}}
\affil[a]{\normalsize Department of Physics, Faculty of Science - UNAM, Mexico}

\author[,b]{F. Brambila-Paz \footnote{Email: fernando.brambila@ciencias.unam.mx; ORCID: 0000-0001-7896-6460}}
\affil[b]{\normalsize Department of Mathematics, Faculty of Science - UNAM, Mexico}

\author[,c]{R. Montufar-Chaveznava
  \footnote{Email: montufar@unam.mx; ORCID: 0000-0001-7630-6207 }}
\affil[c]{\normalsize Department of Geophysics, Faculty of Engineering - UNAM, Mexico}

\date{}

\begin{document}

\maketitle


\begin{abstract}

This paper presents a way to define, classify and accelerate the order of convergence of an uncountable family of fractional fixed point methods, which may be useful to continue expanding the applications of fractional operators. The proposed method to accelerate convergence is used in a fractional iterative method, and with the obtained method are solved simultaneously two nonlinear algebraic systems that depend on time-dependent parameters, and that allow obtaining the temperatures and efficiencies of a hybrid solar receiver. Finally, two uncountable families of fractional fixed point methods are presented, in which the proposed method to accelerate convergence can be implemented.

\end{abstract}

\section{Introduction}

A fractional derivative is an operator that generalizes the ordinary derivative, in the sense that if

\begin{eqnarray*}
\dfrac{d^\alpha}{d x^\alpha},
\end{eqnarray*}

denotes the differential of order $ \alpha $, then $\alpha$ may be considered a parameter, with $\alpha\in \nset{R}$, such that the first derivative corresponds to the particular case $\alpha=1$. On the other hand, a fractional differential equation is an equation that involves at least one differential operator of order $ \alpha $, with $(n-1)< \alpha \leq n$ for some positive integer $n$, and it is said to be a differential equation of order $\alpha$ if this operator is the highest order in the equation. The fractional operators have many representations, but one of their fundamental properties is that they allow retrieving the results of conventional calculus when $\alpha \to n$. So, considering a scalar function $h: \nset{R}^m \to \nset{R}$ and the canonical basis of $\nset{R}^m$ denoted by $\set{\hat{e}_k}_{k\geq 1}$, it is possible to define the following fractional operators of order $\alpha$ using Einstein notation

\begin{eqnarray}
o^\alpha h(x):=\hat{e}_k o_k^\alpha h(x).
\end{eqnarray}

Therefore, denoting by  $\partial_k^n$ the partial derivative of order $n$ applied with respect to the $k$-th component of the vector $x$, using the previous operator it is possible to define the following set of fractional operators

\begin{eqnarray}\label{eq:0}
\Op_{x,\alpha}^n(h):=\set{o^\alpha \ : \ \exists \partial_k^n h(x)  \  \mbox{ and } \ \lim_{\alpha \to n}o_k^\alpha h(x)=\partial_k^n h(x) \ \forall k\geq 1 },
\end{eqnarray}

which may be considered as a generating set of sets of fractional tensor operators. For example, considering $\alpha,n\in \nset{R}^m$ with $\alpha=\hat{e}_k\alpha_k$ and $n=\hat{e}_k n_k$, it is possible to define the following set of fractional tensor operators

\begin{eqnarray}
\Op_{x,\alpha}^{n}(h):=\set{o^\alpha \ : \ o^\alpha \in \Op_{x,\alpha_1}^{n_1}(h) \times \Op_{x,\alpha_2}^{n_2}(h)\times \cdots \times \Op_{x,\alpha_m}^{n_m}(h) }.
\end{eqnarray}

One of the most famous fixed point methods is the well-known Newton-Raphson method. However, it sometimes goes unnoticed that this method has the following problem related to finding roots of polynomials in the complex space: If it is necessary to find a complex root $\xi$ of a polynomial using the Newton-Raphson method, with $\xi\in \nset{C}\setminus\nset{R}$, a complex initial condition $x_0$ must be provided, and if a suitable initial condition is selected, this will lead to a complex solution, but there is also the possibility that this may lead to a real solution. If the root obtained is real, it is necessary to change the initial condition and expect that this will lead to a complex solution, otherwise, it is necessary to change the value of the initial condition again, this process is repeated until it finally converges to a complex solution. The process described above is very similar to what happens when different values $\alpha$ are used in fractional operators until we find a solution that meets some established criterion. Considering the Newton-Raphson method from the perspective of fractional calculus, it is possible to consider that an order $\alpha$ remains fixed, in this case $\alpha=1$, and the initial conditions $x_0$ are varied until found a solution $\xi$ that fulfills an established criterion. It is necessary to mention that considering a relationship between fractional calculus and the Newton-Raphson method may seem somewhat forced at first, but the latter is characterized by the fact that when it generates divergent sequences of complex numbers, it can sometimes lead to the creation of fractals \cite{tatham}, and this feature is complemented quite well with the fact that the orders of the fractional derivatives seem to be closely related to the fractal dimension  \cite{brambila2017fractal}. Based on the above, it is possible to consider inverting the behavior of the order $\alpha=1$ of the derivative  and the initial condition $x_0$, that is, leaving the initial condition $x_0$ fixed and varying the order $\alpha$ of the derivative, thus obtaining the \textbf{fractional Newton-Raphson method}  \cite{torres2021fracnewrap}, which is nothing other than the Newton-Raphson method using any definition of the fractional derivative that fits the function whose zeros want to be determined.

\begin{figure}[!ht]
        \begin{subfigure}[c]{0.35\textwidth}
        \centering
 \includegraphics[width=\textwidth, height=0.6\textwidth]{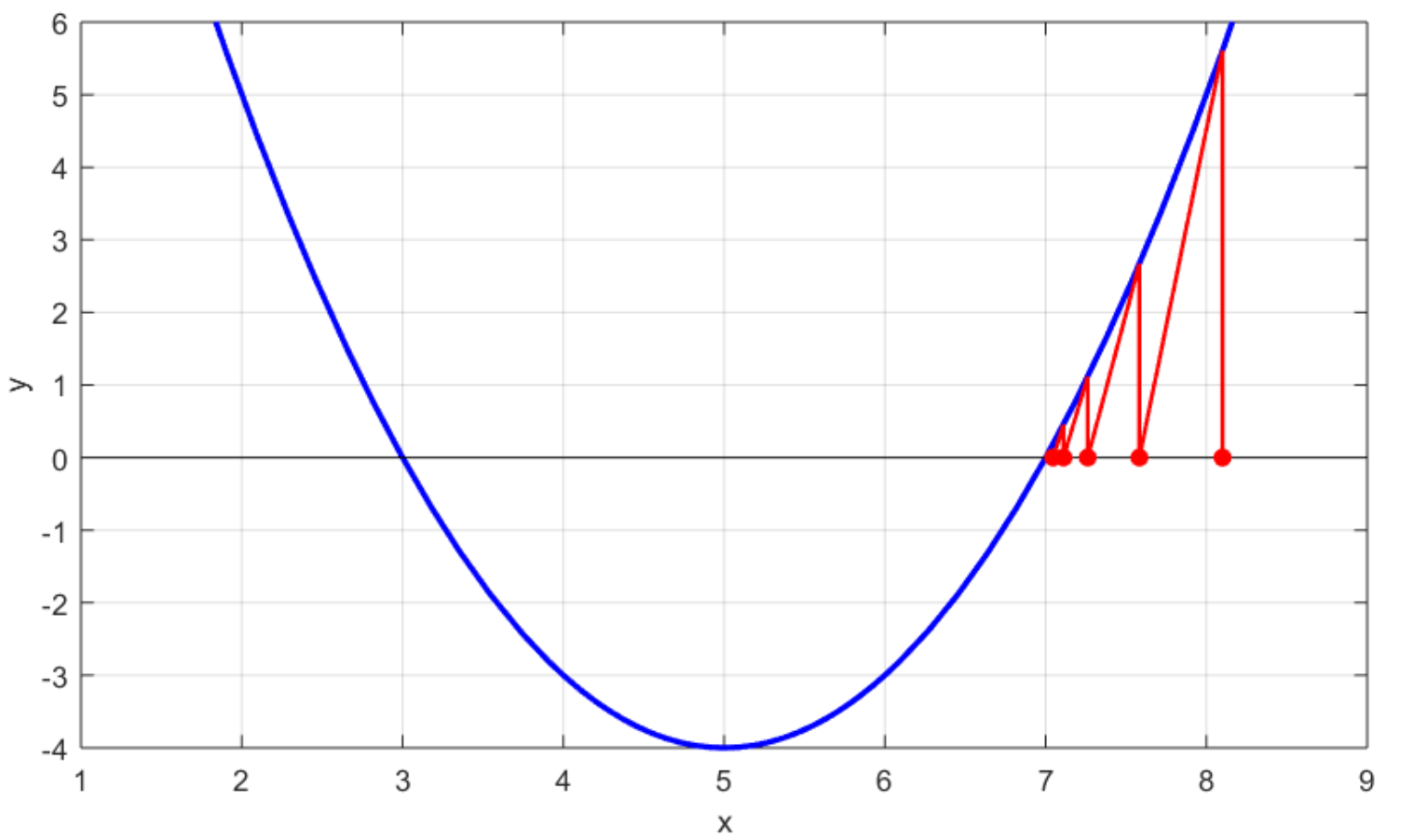}      
    \caption*{a) $\alpha=-0.77$}
    \end{subfigure}
        \begin{subfigure}[c]{0.35\textwidth}
        \centering
 \includegraphics[width=\textwidth, height=0.6\textwidth]{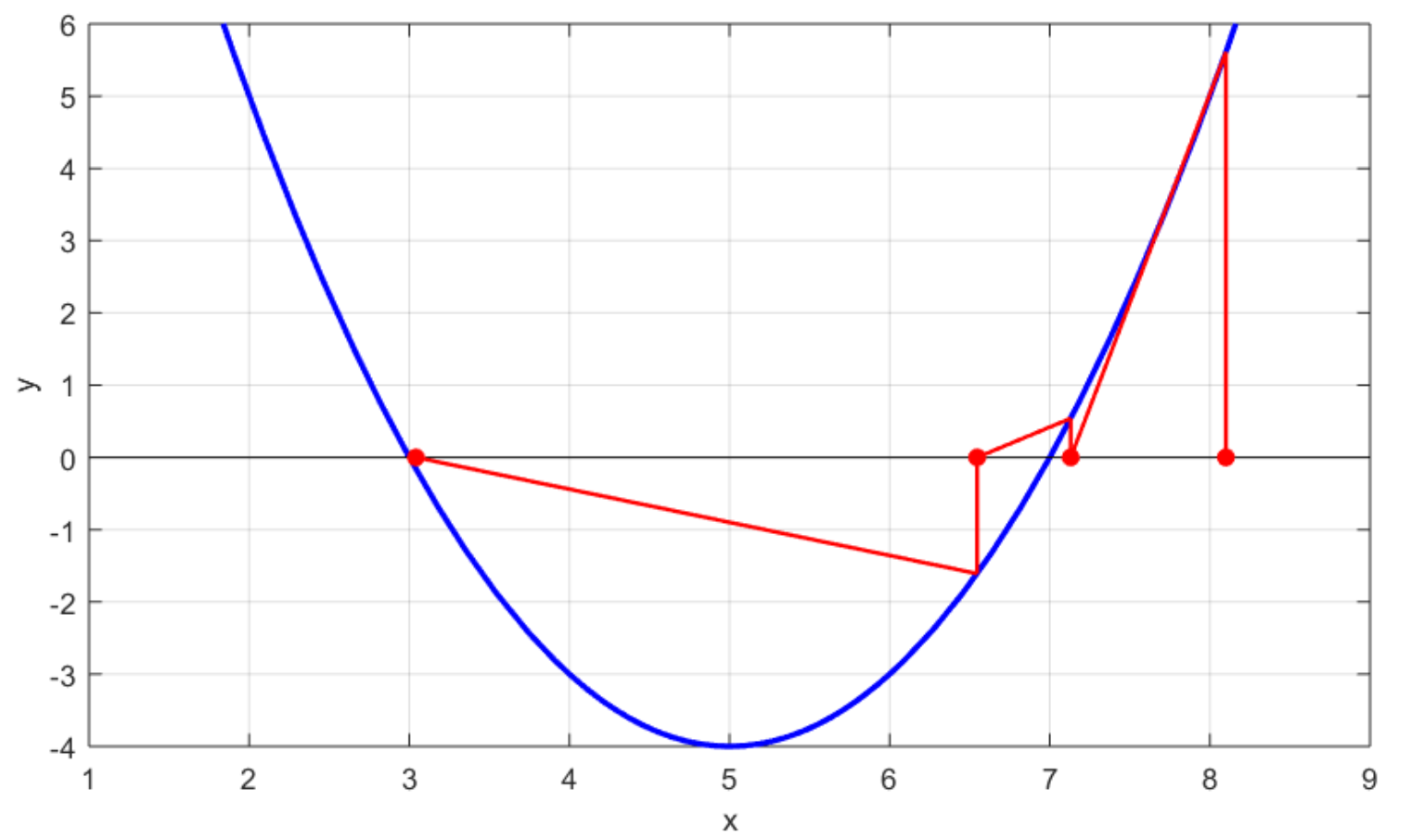}     
    \caption*{b) $\alpha=-0.32$}
    \end{subfigure} 
    \centering
    \begin{subfigure}[c]{0.35\textwidth}
    \centering
 \includegraphics[width=\textwidth, height=0.6\textwidth]{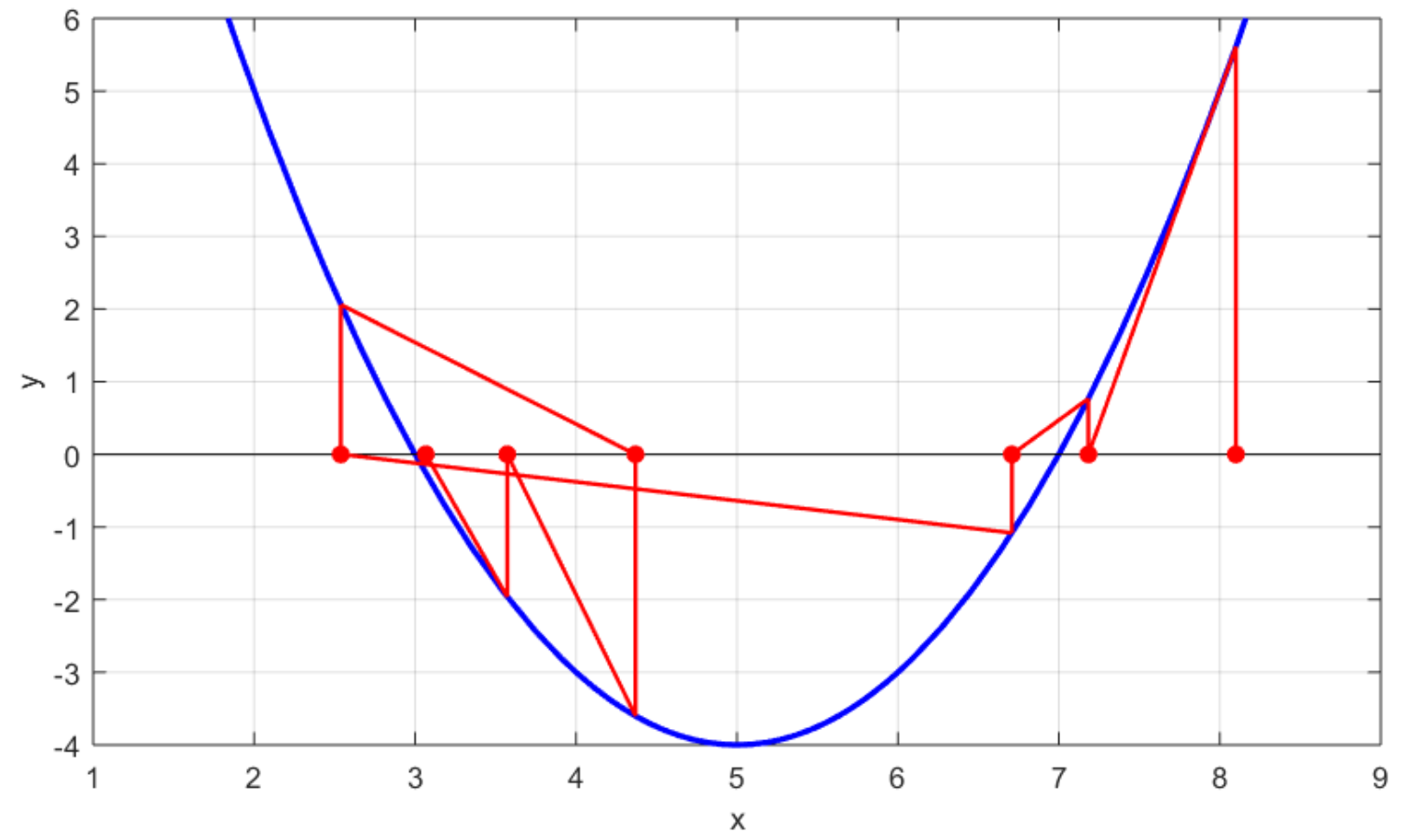}      
    \caption*{c) $\alpha=0.19$}
    \end{subfigure}
    \begin{subfigure}[c]{0.35\textwidth}
    \centering
 \includegraphics[width=\textwidth, height=0.6\textwidth]{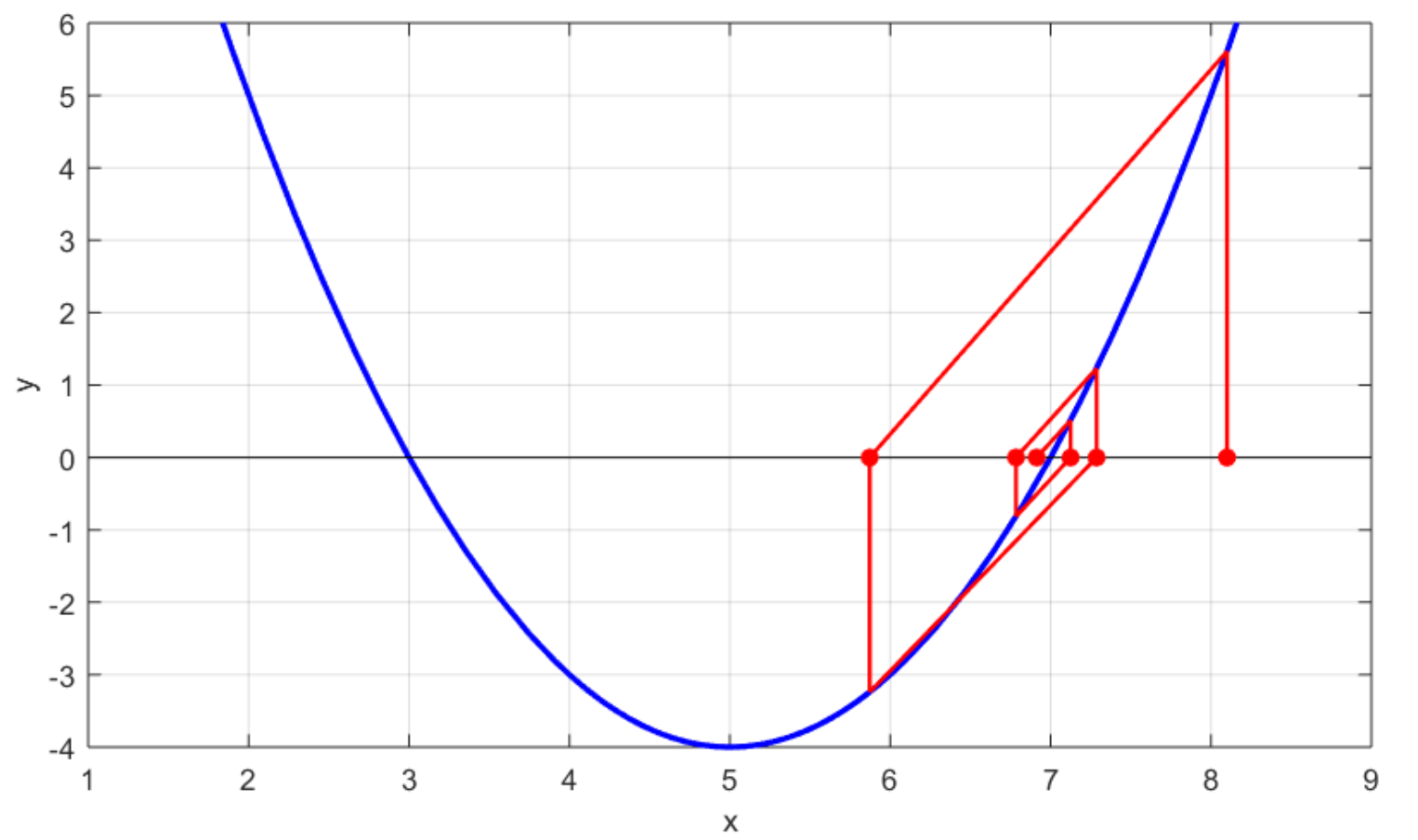}      
    \caption*{d) $\alpha=1.87$}
    \end{subfigure}        
        \caption{Illustrations of some trajectories generated by the fractional Newton-Raphson method for the same initial condition $x_0$ but with different orders $\alpha$ of the fractional operator used \cite{torres2021fracnewrap}.}\label{fig:01}
\end{figure}

Before continuing, it is necessary to mention that due to the large number of fractional operators that can exist, some sets need to be defined to fully characterize the fractional Newton-Raphson method. So, considering a function $h: \Omega \subset \nset{R}^m \to \nset{R}^m$, it is possible to define the following set of fractional operators

\begin{eqnarray}
{}_m\Op_{x,\alpha}^n(h):=\set{ o^\alpha \ : \ o^\alpha\in \Op_{x,\alpha}^n\left([h]_j \right) \ \forall j \mbox{ with } 1\leq j\leq m},
\end{eqnarray}

where $[h]_k: \Omega \subset \nset{R}^m \to \nset{R}$ denotes the $k$-th component of the function $h$. So, it is possible to define the following set of matrices

\begin{eqnarray}
{}_{m}\Ma_{x,\alpha}^{n}(h):=\set{A=A(o^\alpha) \ : \   o^\alpha \in {}_m\Op_{x,\alpha}^n(h) \ \mbox{ and } \ A(x)=\left([A]_{jk}(x) \right):= \left( o_k^\alpha [h]_j(x)\right)  },
\end{eqnarray}

and therefore, the fractional Newton-Raphson method can be defined and classified through the following set of  matrices

\begin{eqnarray}
{}_{m}\IMa_{x,\alpha}^{1}(h):=\set{A \ : \   \exists A^{-1} \in  {}_mM_{x,\alpha}^{1}(h) },
\end{eqnarray}

and as a consequence, if $\Phi_{FNR}$ denotes the iteration function of the fractional Newton-Raphson method, it is possible to obtain the following result:

\begin{eqnarray}
\mbox{Let }\alpha_0\in \nset{R}\setminus \nset{Z} \Rightarrow \forall A \in {}_{m}\IMa_{x,\alpha_0}^{1}(h) \hspace{0.1cm}  \exists\Phi_{FNR}=\Phi_{FNR}(\alpha_0, A) \  \therefore  \ \forall A \hspace{0.1cm}  \exists  \set{\Phi_{FNR}(\alpha, A) \ : \ \alpha \in \nset{R}\setminus \nset{Z}}.  
\end{eqnarray}

The change from leaving the initial condition $x_0$ fixed and varying the order $\alpha$ of the derivative, although seemingly simple, gives the fractional Newton-Raphson method the ability to partially solve the intrinsic problem associated with classical fixed point methods, which is that in general, to find $N$ zeros of a function, $N$ initial conditions must be provided. This is because by varying the order $\alpha$ of the fractional derivative, the fractional Newton-Raphson method can find $N$ zeros of a function using a single initial condition as shown in the Figure \ref{fig:01}. It is necessary to consider that mentioned above is also valid for any fixed point method that implements fractional operators in some way, which may be named as \textbf{fractional fixed point methods} or \textbf{fractional iterative methods}.

To finish this section, it is necessary to mention that the applications of fractional operators have spread to different fields of science such as finance \cite{safdari2015radial,sabatelli2002waiting}, economics \cite{traore2020model,torres2020nonlinear}, number theory through the Riemann zeta function \cite{guariglia2021fractional,torres2021zeta} and in engineering with the study for the manufacture of hybrid solar receivers \cite{de2021fractional,torres2020reduction,torres2020fracpseunew,torres2021codefracpseudo}. It should be mentioned that there is also a growing interest in fractional operators and their properties for the solution of nonlinear systems \cite{wang2021derivative,gdawiec2019visual,akgul2019fractional,torres2021fracnewrap,torres2020fracsome,torres2021fracnewrapaitken}, which is a classic problem in mathematics, physics and engineering, which consists of finding the set of zeros of a function $f:\Omega \subset \nset{R}^n \to \nset{R}^n$, that is,

\begin{eqnarray}\label{eq:1-001}
\set{\xi \in \Omega \ : \ \norm{f(\xi)}=0},
\end{eqnarray}

where $\norm{ \ \cdot \ }: \nset{R}^n \to \nset{R}$ denotes any vector norm, or equivalently

\begin{eqnarray}
\set{\xi \in \Omega \ : \ [f]_k(\xi)=0 \ \forall k\geq 1}.
\end{eqnarray}

Although finding the zeros of a function may seem like a simple problem, it is generally necessary to use numerical methods of the iterative type to solve it. So, considering that fractional iterative methods can find $N$ solutions of a system using a single initial condition, this article shows an alternative way to the Aitken's method to accelerate the order of convergence of a family of fractional fixed point methods, which consists of implementing a function in the order of the fractional operators involved, with which it is possible to obtain an order of convergence at least quadratic.

\section{Fixed Point Method}

Let  $\Phi:\nset{R}^n \to \nset{R}^n$ be a function. It is possible to build a sequence $\set{x_i}_{\geq 1}$  by defining the following iterative method

\begin{eqnarray}\label{eq:2-001}
x_{i+1}:=\Phi(x_i), & i=0,1,2,\cdots,
\end{eqnarray}

if it is fulfilled that $x_i\to \xi\in \nset{R}^n$ and if the function $\Phi$ is continuous around $\xi$, we obtain that

\begin{eqnarray}\label{eq:2-002}
\ds \xi=\lim_{i\to \infty}x_{i+1}=\lim_{i\to \infty}\Phi(x_i)=\Phi\left(\lim_{i\to \infty}x_i \right)=\Phi(\xi),
\end{eqnarray}

the above result is the reason by which the method \eqref{eq:2-001} is known as the \textbf{fixed point method}. Furthermore, the function $\Phi$ is called an \textbf{iteration function}. The following corollary allows characterizing the order of convergence of an iteration function $ \Phi $ through its \textbf{Jacobian matrix} $\Phi^{(1)}$ \cite{torres2021fracnewrapaitken}:

\begin{corollary}\label{cor:2-001}
Let $\Phi:\nset{R}^n \to \nset{R}^n$ be an iteration function. If $\Phi$ defines a sequence $\set{x_i}_{\geq 1}$ such that $x_i\to \xi\in \nset{R}^n$. So, $\Phi$ has an \textbf{order of convergence}  of order (at least) $p$  in $B(\xi;\delta)$, where

\begin{eqnarray}\label{eq:2-003}
p:=\left\{
\begin{array}{cc}
1 , &\ds \mbox{If } \lim_{x\to \xi}\norm{\Phi^{(1)}(x)}\neq 0  \vspace{0.2cm}\\
2, &\ds \mbox{If } \lim_{x\to \xi}\norm{\Phi^{(1)}(x)}= 0  \\
\end{array}\right. .
\end{eqnarray}

\end{corollary}

\section{Riemann-Liouville Fractional Operators}

One of the fundamental operators of fractional calculus is the operator \textbf{Riemann-Liouville fractional integral},  which is defined as follows \cite{hilfer00,oldham74}

\begin{eqnarray}
\ifr{}{a}{I}{x}{\alpha}f(x):=\dfrac{1}{\gam{\alpha}}\int_a^x (x-t)^{\alpha-1}f(t)dt,
\end{eqnarray}

which is a fundamental piece to build the operator \textbf{Riemann-Liouville fractional derivative},  which is defined as follows \cite{hilfer00,kilbas2006theory}

\begin{eqnarray}\label{eq:3-001}
\normalsize
\begin{array}{c}
\ifr{}{a}{D}{x}{\alpha}f(x) := \left\{
\begin{array}{cc}
\ds \ifr{}{a}{I}{x}{-\alpha}f(x), &\mbox{if }\alpha<0 \vspace{0.1cm}\\  
\ds \dfrac{d^n}{dx^n}\left( \ifr{}{a}{I}{x}{n-\alpha}f(x)\right), & \mbox{if }\alpha\geq 0
\end{array}
\right.
\end{array}, 
\end{eqnarray}

where  $ n = \lceil \alpha \rceil$ and $\ifr{}{a}{I}{x}{0}f(x):=f(x)$. Applying the operator \eqref{eq:3-001} with $a=0$ to the  function $ x^{\mu} $, with $\mu> -1$, we obtain the following result \cite{torres2021fracnewrapaitken}:

\begin{eqnarray}\label{eq:3-002}
\ifr{}{0}{D}{x}{\alpha}x^\mu = 
 \dfrac{\gam{\mu+1}}{\gam{\mu-\alpha+1}}x^{\mu-\alpha}, & \alpha\in \nset{R}\setminus \nset{Z}.
\end{eqnarray}

\section{Fractional Fixed Point Method}

Let $f:\Omega \subset \nset{R}^n \to \nset{R}^n$ be a function with a point $\xi\in \Omega$ such that $\norm{f(\xi)}=0$ . So, considering an iteration function $\Phi:(\nset{R}\setminus \nset{Z})\times \nset{R}^n\to \nset{R}^n$, the iteration function of a fractional iterative method may be written in general form as follows

\begin{eqnarray}\label{eq:4-001}
\Phi(\alpha,x):=x-A_{g,\alpha}(x)f(x),& \alpha\in \nset{R}\setminus \nset{Z},
\end{eqnarray}

where $A_{g,\alpha}$ is a matrix that depends, in at least one of its entries, on fractional operators of order $\alpha$ applied to some function $g:\nset{R}^n \to \nset{R}^n$, whose particular case occurs when $g=f$. So, it is possible to define in a general way a fractional fixed point method as follows

\begin{eqnarray}\label{eq:4-002}
x_{i+1}:=\Phi(\alpha,x_i), & i=0,1,2,\cdots.
\end{eqnarray}

If it is fulfilled that $x_i\to \xi\in \Omega$, it is possible to define the following set

\begin{eqnarray}
\Conv(\xi):=\set{\Phi \ : \ \lim_{x\to \xi}\Phi(\alpha,x)=\xi},
\end{eqnarray}

which may be interpreted as the set of fractional fixed point methods that define a convergent sequence $\set{x_i}_{i\geq 1}$ to the value $\xi \in \Omega$. Considering the \mref{Corollary}{\ref{cor:2-001}}, as well as the \textbf{Proposition 1} from the reference \cite{torres2021fracnewrapaitken}, it is possible to define the following sets to classify the order of convergence of some fractional iterative methods

\begin{eqnarray}
\Ord^1(\xi):=\set{\Phi \in \Conv(\xi) \ : \ \lim_{x \to \xi}\norm{\Phi^{(1)}(\alpha,x)}\neq 0 },
\end{eqnarray}

\begin{eqnarray}
\Ord^2(\xi):=\set{\Phi \in \Conv(\xi) \ : \ \lim_{x \to \xi}\norm{\Phi^{(1)}(\alpha,x)}= 0 },
\end{eqnarray}

\begin{eqnarray}
\ord^1(\xi):=\set{\Phi \in \Conv(\xi) \ : \ \lim_{x \to \xi}A_{g,\alpha}(x)\neq \left(f^{(1)}(\xi)\right)^{-1} \ \mbox{ or } \ \lim_{\alpha \to 1} A_{g,\alpha}(\xi)\neq \left(f^{(1)}(\xi)\right)^{-1} },
\end{eqnarray}

\begin{eqnarray}
\ord^2(\xi):=\set{\Phi \in \Conv(\xi) \ : \ \lim_{x \to \xi}A_{g,\alpha}(x)= \left(f^{(1)}(\xi)\right)^{-1} \ \mbox{ or } \ \lim_{\alpha \to 1} A_{g,\alpha}(\xi)= \left(f^{(1)}(\xi)\right)^{-1}}.
\end{eqnarray}

On the other hand, considering that depending on the nature of the function $f$ there exist cases in which the Newton-Raphson method can present an order of convergence (at least) linear \cite{torres2021fracnewrapaitken}, it is possible to obtain the following relations between the previous sets

\begin{eqnarray}
\ord^1(\xi) \subset \Ord^1(\xi)& \mbox{ and } &\ord^2(\xi)\subset \Ord^1(\xi)\cup \Ord^2(\xi),
\end{eqnarray}

with which it is possible to define the following sets

\begin{eqnarray}
\Ord_2^1(\xi) := \ord^2(\xi)   \cap  \Ord^1(\xi) & \mbox{ and } & \Ord_2^2(\xi):= \ord^2(\xi) \cap \Ord^2(\xi). 
\end{eqnarray}

\subsection{Acceleration in the Order of Convergence of the Set $\Ord_2^1(\xi)$}

Let $f:\Omega \subset \nset{R}^n \to \nset{R}^n$ be a function with a point $\xi\in \Omega$ such that $\norm{f(\xi)}=0$, and denoting by $\Phi_{NR}$ to the iteration function of the Newton-Raphson method, it is possible to define the following set of functions

\begin{eqnarray}
\Ord_{NR}^2(\xi):=\set{ f \ : \  \lim_{x\to \xi}\norm{\Phi_{NR}^{(1)}(x)}=0 }.
\end{eqnarray}

So, it is possible to define the following corollary:

\begin{corollary}\label{cor:4-001}
Let $f:\Omega \subset \nset{R}^n \to \nset{R}^n$ be a function with a point $\xi\in \Omega$ such that $f\in \Ord_{NR}^2(\xi)$, and let $\Phi$ be a iteration function given by the equations \eqref{eq:4-001} such that $\Phi\in \Ord_2^1(\xi)$. So, it is possible to replace the order $\alpha$ of the fractional operators of the matrix $A_{g,\alpha}$ by the following function

\begin{eqnarray}\label{eq:4-005}
\alpha_f([x]_k,x):=\left\{
\begin{array}{cc}
\alpha ,& \mbox{if \hspace{0.1cm}}\abs{[x]_k}\neq 0 \mbox{ \hspace{0.1cm}and\hspace{0.1cm} } \norm{f(x)}> \delta \vspace{0.1cm}\\
1,& \mbox{if \hspace{0.1cm}}\abs{[x]_k} = 0 \mbox{ \hspace{0.1cm}or\hspace{0.1cm} } \norm{f(x)}\leq \delta
\end{array}
\right. ,
\end{eqnarray}

obtaining a new matrix that may be denoted as follows

\begin{eqnarray}\label{eq:4-006}
A_{g,\alpha_f}(x)=\left([A_{g,\alpha_f}]_{jk}(x) \right), & \alpha\in \nset{R}\setminus \nset{Z},
\end{eqnarray}

and that guarantees that there exists a set $\Omega_\xi \subset B(\xi;\delta)$ such that $\Phi\in \Ord_2^2(\xi)$ in $\Omega_\xi$.

\end{corollary}

It is necessary to mention that the origin of the function \eqref{eq:4-005} arises from the need to accelerate the order of convergence of the fractional Newton-Raphson method, which generated the method known as the \textbf{fractional Newton method}, whose matrix $A_{g,\alpha_f}$ corresponds to a particular case in which $g=f$ \cite{torres2020fracsome,torres2021fracnewrap,torres2021fracnewrapaitken}. Finally, for practical purposes, it may be defined that if a fractional iterative method $\Phi \in \Ord_2^1(\xi)$ uses the function \eqref{eq:4-005}, it may be called a \textbf{fractional iterative method accelerated}.

\section{Equations of a Hybrid Solar Receiver}

Considering the notation

\begin{eqnarray*}
s=(T_{cell},T_{hot},T_{cold},\eta_{cell},\eta_{TEG})^T:=\left([x]_1,[x]_2,[x]_3,[x]_4,[x]_5 \right)^T,
\end{eqnarray*}

the following expressions

\begin{eqnarray*}
\begin{array}{c}
\begin{array}{lll}
a_0=\dfrac{2*r_{intercon}}{\cdot \sqrt{f^*\cdot A_{TEG}}\left(b\cdot \sqrt{f^*}+\sqrt{A_{TEG}} \right) }, & a_1=\eta_{opt}\cdot C_g \cdot DNI, & a_2=r_{cell}+r_{sol}+A_{cell}\left(\dfrac{r_{cop}+r_{cer}}{A_{TEG}}+a_0 \right)
\end{array}\vspace{0.2cm}\\
\begin{array}{llll}
a_3=\dfrac{A_{cell}\cdot l}{f^*\cdot A_{TEG}\cdot k_{TEG}}, & a_4=T_{air}, &a_5=A_{cell}\left( \dfrac{r_{cer}}{A_{TEG}}+R_{heat\_exch}+a_0 \right),&a_6=-\eta_{cell,ref}\cdot \gamma_{cell}
\end{array}\vspace{0.2cm}\\
\begin{array}{lll}
 a_7=\eta_{cell,ref}\left(1+25 \cdot \gamma_{cell} \right), &a_8=\sqrt{1+ZT}, &a_9=273.15
\end{array}
\end{array},
\end{eqnarray*}

and the following particular values \cite{rodrigo2019performance}

\begin{eqnarray*}
\left\{
\begin{array}{lll}
\eta_{opt}=0.85, & r_{intercon}=2.331\times 10^{-7} ,&C_g=800\\
    A_{cell}=9\times 10^{-6}     , &  R_{heat\_exch}=0.5   , & A_{TEG}=5.04 \times 10^{-5}     \\
 \eta_{cell,ref}=0.43,& r_{cell}=3\times 10^{-6}   , &    f^*=0.7  \\
 \gamma_{cell}=4.6\times 10^{-4}, &  r_{sol}=1.603\times 10^{-6}  , &  b=5\times 10^{-4}    \\
r_{cop}=7.5\times 10^{-7}, & r_{cer}=8 \times 10^{-6},&    l=5\times 10^{-4}   \\
k_{TEG}=1.5 ,&   ZT=1   & 
\end{array}
\right..
\end{eqnarray*}

it is possible to define the following system of equations that corresponds to the combination of a solar photovoltaic system with a thermoelectric generator system \cite{bjork2015performance,bjork2018maximum}, which is named as a \textbf{hybrid solar receiver}

\begin{eqnarray}\label{eq:5-001}
\left\{
\begin{array}{l}
\left[x\right]_1=[x]_2+a_1\cdot a_2\left( 1-[x]_4 \right)\\
\left[x\right]_2=[x]_3+a_1\cdot a_3 \left( 1-[x]_4\right)\left(1-[x]_5\right)\\
\left[x\right]_3=a_4+a_1\cdot a_5 \left( 1-[x]_4\right)\left(1-[x]_5\right)\\
\left[x\right]_4 =a_6[x]_1+a_7\\
\left[x\right]_5=(a_8-1)\left(1-\dfrac{[x]_3+a_9}{[x]_2+a_9} \right)\left(a_8+ \dfrac{[x]_3+a_9}{[x]_2+a_9}\right)^{-1}
\end{array}\right.,
\end{eqnarray}

whose deduction, as well as details about its interpretation, may be found in the reference \cite{rodrigo2019performance}. Using the system of equations \eqref{eq:5-001}, it is possible to define a  function $f_1:\Omega \subset \nset{R}^5\to \nset{R}^5$, that is,

\begin{eqnarray}\label{eq:5-002}
f_1(s):=\begin{pmatrix}
\left[x\right]_1-[x]_2-a_1\cdot a_2\left( 1-[x]_4 \right)\\
\left[x\right]_2-[x]_3-a_1\cdot a_3 \left( 1-[x]_4\right)\left(1-[x]_5\right)\\
\left[x\right]_3-a_4-a_1\cdot a_5 \left( 1-[x]_4\right)\left(1-[x]_5\right)\\
\left[x\right]_4 -a_6[x]_1-a_7\\
\left[x\right]_5-(a_8-1)\left(1-\dfrac{[x]_3+a_9}{[x]_2+a_9} \right)\left(a_8+ \dfrac{[x]_3+a_9}{[x]_2+a_9}\right)^{-1}
\end{pmatrix}, 
\end{eqnarray}

which depends on two parameters, the direct normal irradiance ($DNI$) and the ambient temperature ($T_{air}$). These parameters are measured in real-time at certain times of the day \cite{rodrigo2019performance}, and it is necessary to calculate a new solution of the system \eqref{eq:5-001} for each new pair of parameters, that is,

\begin{eqnarray*}
(DNI,T_{air})\overset{f_1}{\longrightarrow} s \in \nset{R}^5.
\end{eqnarray*}

However, to simplify the task of finding the solutions of the function \eqref{eq:5-002}, it is possible through the consecutive substitution of the variables $[x]_1, \ [x]_4, \ [x]_5$ and some algebraic simplifications, to obtain the following transcendental system \cite{torres2020reduction}

\begin{eqnarray}\label{eq:5-003}
\left\{
\begin{array}{l}
\left[x\right]_2=[x]_3-a_1\cdot a_3 \dfrac{\left( a_6 [x]_2 + a_7 - 1 \right) \left( a_8 \left([x]_3+ a_9 \right) + \left([x]_2+ a_9\right)  \right)  }{(1+a_1 a_2 a_6 ) \left( a_8 \left( [x]_2 + a_9  \right) + \left([x]_3+ a_9\right) \right) } \vspace{0.1cm} \\
\left[x \right]_3=a_4-a_1\cdot a_5 \dfrac{\left( a_6 [x]_2 + a_7 - 1 \right) \left( a_8 \left([x]_3+ a_9 \right) + \left([x]_2+ a_9\right)  \right)  }{(1+a_1 a_2 a_6 ) \left( a_8 \left( [x]_2 + a_9  \right) + \left([x]_3+ a_9\right) \right) } 
\end{array}\right.,
\end{eqnarray}

whose solution allows to know the values of the variables $[x]_1,[x]_4$ and $[x]_5$ through the following equations

\begin{eqnarray}\label{eq:5-004}
\left\{
\begin{array}{l}
\left[x\right]_1=\dfrac{[x]_2 - a_1 a_2 (a_7 - 1)}{1+a_1 a_2 a_6} \vspace{0.1cm}\\
\left[ x \right]_4=\dfrac{a_6 \left( a_1 a_2 + [x]_2\right) + a_7}{1+ a_1 a_2 a_6 } \vspace{0.1cm}\\
\left[ x \right]_5=\dfrac{(a_8 - 1)\left( [x]_2 - [x]_3 \right)}{a_8 \left([x]_2+a_9\right) + \left( [x]_3+a_9 \right)}
\end{array}\right..
\end{eqnarray}

Using the system of equations \eqref{eq:5-003}, it is possible to define a  function $f_2:\Omega \subset \nset{R}^2\to \nset{R}^2$, that is,

\begin{eqnarray}\label{eq:5-005}
f_2(x):=\begin{pmatrix}
\left[x\right]_2-[x]_3+a_1\cdot a_3 \dfrac{\left( a_6 [x]_2 + a_7 - 1 \right) \left( a_8 \left([x]_3+ a_9 \right) + \left([x]_2+ a_9\right)  \right)  }{(1+a_1 a_2 a_6 ) \left( a_8 \left( [x]_2 + a_9  \right) + \left([x]_3+ a_9\right) \right) } \vspace{0.1cm} \\
\left[x \right]_3-a_4+a_1\cdot a_5 \dfrac{\left( a_6 [x]_2 + a_7 - 1 \right) \left( a_8 \left([x]_3+ a_9 \right) + \left([x]_2+ a_9\right)  \right)  }{(1+a_1 a_2 a_6 ) \left( a_8 \left( [x]_2 + a_9  \right) + \left([x]_3+ a_9\right) \right) } 
\end{pmatrix}, 
\end{eqnarray}

and then finding the solutions of the function \eqref{eq:5-005}, through the equations \eqref{eq:5-004}, it is possible to construct the solutions of the function \eqref{eq:5-002}.

\subsection{Solutions of the Equations of a Hybrid Solar Receiver}

To solve the equation \eqref{eq:5-005} and at the same time solve the equation \eqref{eq:5-002}, a fractional fixed point method will be used as well as its accelerated version through the function \eqref{eq:4-005}. Before continuing, it is necessary to mention that for some definitions of fractional operators it is fulfilled that the derivative of order $\alpha$ of a constant is different from zero (for example: Riesz, Grünwald–Letnikov, Riemann-Liouville, etc.\cite{miller93,hilfer00,oldham74,kilbas2006theory,brambila2017fractal}), that is,

\begin{eqnarray}\label{eq:5-006}
\partial_k^\alpha c :=\der{\partial}{[x]_k}{\alpha}c \neq 0 , & c=constant.
\end{eqnarray}

So, considering a function $f:\Omega \subset \nset{R}^n \to \nset{R}^n$ with a point $\xi \in \Omega$ such that $\norm{f(\xi)}=0$, the Riemann-Liouville fractional derivative given by the equation \eqref{eq:3-002}, and an iteration function $\Phi:(\nset{R}\setminus \nset{Z})\times \nset{R}^n \to \nset{R}^n$, it is possible to define the following fractional fixed point method

\begin{eqnarray}\label{eq:5-007}
x_{i+1}:=\Phi(\alpha,x_i)=x_i-A_{g_f,\beta}(x_i)f(x_i), & i=0,1,2,\cdots,
\end{eqnarray}

where $ A_{g_f,\beta}(x_i) $ is given by the following expression

\begin{eqnarray}\label{eq:5-008}
A_{g_f,\beta}(x_i)=\left([A_{g_f,\beta}]_{jk}(x_i) \right) :=\left( \partial_k^{\beta(\alpha,[x]_k)}[g_f]_{j}(x)  \right)_{x_i}^{-1}, & \alpha\in \nset{R}\setminus \nset{Z},
\end{eqnarray}

with $ g_{f}(x) $  and $ \beta (\alpha, [x]_k) $ functions defined as follows

\begin{eqnarray}
g_{f}(x):=f(x_i)+f^{(1)}(x_i)x& \mbox{ and }
&
\beta(\alpha,[x]_k):=\left\{
\begin{array}{cc}
\alpha, &\mbox{if \hspace{0.1cm} } \abs{  [x]_k }\neq 0 \vspace{0.1cm}\\
1,& \mbox{if \hspace{0.1cm}  }  \abs{  [x]_k  }=0
\end{array}\right..
\end{eqnarray}

The fractional iterative method given by the equation \eqref{eq:5-007} is named the \textbf{fractional quasi-Newton method} \cite{torres2020fracsome}, if it is assumed that $\Phi \in \Conv(\xi)$ then $\Phi \in \ord^1(\xi)$. Furthermore, the method fulfills the following condition

\begin{eqnarray}
\lim_{\alpha \to 1}\partial_k^{\beta(\alpha,[x_i]_k)} [g]_j(x_i)=\partial_k[f]_j(x_i), & 1\leq j,k\leq n,
\end{eqnarray}

and as a consequence $\Phi \in \Ord_{2}^1(\xi)$. So, if it is assumed that $f\in \Ord_{NR}^2(\xi)$, by the \mref{Corollary}{\ref{cor:4-001}} it is possible to construct the \textbf{fractional quasi-Newton method accelerated} using the following matrix

\begin{eqnarray}\label{eq:5-009}
A_{g_f,\alpha_f}(x_i)=\left([A_{g_f,\alpha_f}]_{jk}(x_i) \right) :=\left( \partial_k^{\alpha_f([x]_k,x)}[g_f]_{j}(x)  \right)_{x_i}^{-1}, & \alpha\in \nset{R}\setminus \nset{Z}.
\end{eqnarray}

Before continuing, it is necessary to mention that a description of the algorithm that must be implemented when working with a fractional iterative method given by the equation \eqref{eq:4-002} may be found in the reference \cite{torres2020fracsome}. On the other hand, simplified examples of how the methods given by the matrices \eqref{eq:5-008} and \eqref{eq:5-009} should be programmed may be found in the references \cite{torres2021codefracquasi,torres2021codeaccelfracquasi}. Using the fractional fixed point methods defined by the matrices \eqref{eq:5-008} and \eqref{eq:5-009}, we proceed to find three solutions of the function \eqref{eq:5-005}  leaving the following fixed values

\begin{eqnarray*}
\delta=13& \mbox{ and } & x_0=(3000,3000)^T.
\end{eqnarray*}

\begin{example}

Considering by hypothesis that $f_2\in \Ord_{NR}^2(\xi)$, and using the following values

\begin{eqnarray*}
DNI=900, &
Tair=20, & \alpha=0.89825,
\end{eqnarray*}

the following iterations are obtained by using the fractional iterative methods given by the matrices  \eqref{eq:5-008} and \eqref{eq:5-009}.

\begin{itemize}
\item[i)] $A_{g_f,\beta} \Rightarrow \Phi \in \Ord_2^1(\xi)$.

\begin{footnotesize}
\centering
\begin{longtable}{c|cccc|cccccccc}
\toprule
$i$& $[ x_i ]_2$ &$[ x_i ]_3$ &$\norm{x_i-x_{i-1}}_2$&$\norm{f_2(x_i)}_2$&$[ x_i ]_1$&$[ x_i ]_4$&$[ x_i ]_5$& $\norm{f_1(s_i)}_2$ \\
\midrule
    1     & 2048.526273 & 2036.688326 & 1.35E+03 & 2.01E+03 & 2052.245932 & 0.02901075 & 0.00087668 & 2.01E+03 \\
    2     & 1378.380727 & 1357.837031 & 9.54E+02 & 1.33E+03 & 1381.592211 & 0.16166606 & 0.00214528 & 1.33E+03 \\
    3     & 914.5756647 & 887.7554749 & 6.60E+02 & 8.65E+02 & 917.4354426 & 0.25347627 & 0.00391089 & 8.65E+02 \\
    4     & 599.7868499 & 568.5654338 & 4.48E+02 & 5.46E+02 & 602.4079218 & 0.31578871 & 0.00622874 & 5.46E+02 \\
    5     & 390.7721777 & 356.5990844 & 2.98E+02 & 3.34E+02 & 393.2347526 & 0.35716317 & 0.0090235 & 3.34E+02 \\
    6     & 255.3927888 & 219.4044444 & 1.93E+02 & 1.97E+02 & 257.7527048 & 0.38396151 & 0.0120214 & 1.97E+02 \\
    7     & 170.1536777 & 133.2761535 & 1.21E+02 & 1.11E+02 & 172.4489564 & 0.4008346 & 0.01478215 & 1.11E+02 \\
    8     & 118.188164 & 81.23045449 & 7.35E+01 & 5.95E+01 & 120.4440369 & 0.41112117 & 0.01686287 & 5.95E+01 \\
    9     & 87.62585188 & 51.33933793 & 4.27E+01 & 2.99E+01 & 89.85854925 & 0.41717098 & 0.01800683 & 2.99E+01 \\
    10    & 70.31181026 & 35.35034092 & 2.36E+01 & 1.43E+01 & 72.5313783 & 0.42059829 & 0.01823343 & 1.43E+01 \\
    11    & 60.85889363 & 27.58761689 & 1.22E+01 & 6.66E+00 & 63.07129347 & 0.4224695 & 0.01782623 & 6.66E+00 \\
    12    & 55.92035933 & 24.22121073 & 5.98E+00 & 3.07E+00 & 58.12901425 & 0.42344708 & 0.01721438 & 3.07E+00 \\
    13    & 53.49709311 & 22.89305436 & 2.76E+00 & 1.38E+00 & 55.70391046 & 0.42392677 & 0.01672394 & 1.38E+00 \\
    14    & 52.38726485 & 22.39252245 & 1.22E+00 & 6.02E-01 & 54.59324061 & 0.42414646 & 0.01643587 & 6.02E-01 \\
    15    & 51.90534374 & 22.20463447 & 5.17E-01 & 2.55E-01 & 54.11095406 & 0.42424185 & 0.01629349 & 2.55E-01 \\
    16    & 51.70286627 & 22.13313077 & 2.15E-01 & 1.06E-01 & 53.90832305 & 0.42428193 & 0.01622933 & 1.06E-01 \\
    17    & 51.61937072 & 22.10548244 & 8.80E-02 & 4.32E-02 & 53.82476418 & 0.42429846 & 0.01620181 & 4.32E-02 \\
    18    & 51.58529753 & 22.09465371 & 3.58E-02 & 1.76E-02 & 53.79066515 & 0.42430521 & 0.01619031 & 1.76E-02 \\
    19    & 51.5714752 & 22.09037372 & 1.45E-02 & 7.12E-03 & 53.77683234 & 0.42430794 & 0.01618559 & 7.12E-03 \\
    20    & 51.56588734 & 22.0886717 & 5.84E-03 & 2.87E-03 & 53.77124024 & 0.42430905 & 0.01618366 & 2.87E-03 \\
    21    & 51.56363304 & 22.08799215 & 2.35E-03 & 1.16E-03 & 53.76898424 & 0.42430949 & 0.01618288 & 1.16E-03 \\
    22    & 51.56272473 & 22.08772013 & 9.48E-04 & 4.67E-04 & 53.76807524 & 0.42430967 & 0.01618256 & 4.67E-04 \\
    23    & 51.56235904 & 22.08761106 & 3.82E-04 & 1.88E-04 & 53.76770927 & 0.42430975 & 0.01618243 & 1.88E-04 \\
    24    & 51.56221188 & 22.08756728 & 1.54E-04 & 7.56E-05 & 53.767562 & 0.42430978 & 0.01618238 & 7.58E-05 \\
    25    & 51.56215268 & 22.0875497 & 6.18E-05 & 3.04E-05 & 53.76750275 & 0.42430979 & 0.01618236 & 3.05E-05 \\
    26    & 51.56212886 & 22.08754263 & 2.48E-05 & 1.22E-05 & 53.76747891 & 0.42430979 & 0.01618235 & 1.20E-05 \\
    27    & 51.56211928 & 22.08753979 & 9.99E-06 & 4.92E-06 & 53.76746933 & 0.42430979 & 0.01618235 & 4.72E-06 \\
 \bottomrule
\caption{Iterations generated by the fractional quasi-Newton method.}
\end{longtable}
\end{footnotesize}

\item[ii)] $A_{g_f,\alpha_f} \Rightarrow \Phi \in \Ord_2^2(\xi)$ in $\Omega_\xi \subset B(\xi;\delta)$.

\begin{footnotesize}
\centering
\begin{longtable}{c|cccc|cccccccc}
\toprule
$i$& $[ x_i ]_2$ &$[ x_i ]_3$ &$\norm{x_i-x_{i-1}}_2$&$\norm{f_2(x_i)}_2$&$[ x_i ]_1$&$[ x_i ]_4$&$[ x_i ]_5$& $\norm{f_1(s_i)}_2$ \\
\midrule
    1     & 2048.526273 & 2036.688326 & 1.35E+03 & 2.01E+03 & 2052.245932 & 0.02901075 & 0.00087668 & 2.01E+03 \\
    2     & 1378.380727 & 1357.837031 & 9.54E+02 & 1.34E+03 & 1381.592211 & 0.16166606 & 0.00214528 & 1.34E+03 \\
    3     & 914.5756647 & 887.7554749 & 6.60E+02 & 8.65E+02 & 917.4354426 & 0.25347627 & 0.00391089 & 8.65E+02 \\
    4     & 599.7868499 & 568.5654338 & 4.48E+02 & 5.46E+02 & 602.4079218 & 0.31578871 & 0.00622874 & 5.46E+02 \\
    5     & 390.7721777 & 356.5990844 & 2.98E+02 & 3.34E+02 & 393.2347526 & 0.35716317 & 0.0090235 & 3.34E+02 \\
    6     & 255.3927888 & 219.4044444 & 1.93E+02 & 1.97E+02 & 257.7527048 & 0.38396151 & 0.0120214 & 1.97E+02 \\
    7     & 170.1536777 & 133.2761535 & 1.21E+02 & 1.11E+02 & 172.4489564 & 0.4008346 & 0.01478215 & 1.11E+02 \\
    8     & 118.188164 & 81.23045449 & 7.36E+01 & 5.95E+01 & 120.4440369 & 0.41112117 & 0.01686287 & 5.95E+01 \\
    9     & 87.62585188 & 51.33933793 & 4.28E+01 & 2.99E+01 & 89.85854925 & 0.41717098 & 0.01800683 & 2.99E+01 \\
    10    & 70.31181026 & 35.35034092 & 2.36E+01 & 1.43E+01 & 72.5313783 & 0.42059829 & 0.01823343 & 1.43E+01 \\
    11    & 60.85889363 & 27.58761689 & 1.22E+01 & 6.66E+00 & 63.07129347 & 0.4224695 & 0.01782623 & 6.66E+00 \\
    12    & 51.56100988 & 22.08746493 & 1.08E+01 & 1.04E-03 & 53.76635909 & 0.42431001 & 0.01618182 & 1.04E-03 \\
    13    & 51.56211284 & 22.08753788 & 1.11E-03 & 4.13E-09 & 53.76746288 & 0.4243098 & 0.01618235 & 3.03E-07 \\
    \bottomrule
\caption{Iterations generated by the fractional quasi-Newton method accelerated.}
\end{longtable}
\end{footnotesize}

\end{itemize}

\end{example}

\begin{example}

Considering by hypothesis that $f_2\in \Ord_{NR}^2(\xi)$, and using the following values

\begin{eqnarray*}
DNI=574.319,&
Tair=16.832,& \alpha=0.8996,
\end{eqnarray*}

the following iterations are obtained by using the fractional iterative methods given by the matrices  \eqref{eq:5-008} and \eqref{eq:5-009}.

\begin{itemize}
\item[i)] $A_{g_f,\beta} \Rightarrow \Phi \in \Ord_2^1(\xi)$.

\begin{footnotesize}
\centering
\begin{longtable}{c|cccc|cccccccc}
\toprule
$i$& $[ x_i ]_2$ &$[ x_i ]_3$ &$\norm{x_i-x_{i-1}}_2$&$\norm{f_2(x_i)}_2$&$[ x_i ]_1$&$[ x_i ]_4$&$[ x_i ]_5$& $\norm{f_1(s_i)}_2$ \\
\midrule
    1     & 2029.854772 & 2022.247443 & 1.38E+03 & 2.00E+03 & 2032.218723 & 0.03297214 & 0.00056752 & 2.00E+03 \\
    2     & 1351.035349 & 1337.861649 & 9.64E+02 & 1.32E+03 & 1353.07091 & 0.16730757 & 0.00139631 & 1.32E+03 \\
    3     & 884.5286725 & 867.3584839 & 6.63E+02 & 8.49E+02 & 886.3385526 & 0.25962723 & 0.00256042 & 8.49E+02 \\
    4     & 570.3098992 & 550.3428213 & 4.46E+02 & 5.32E+02 & 571.9677708 & 0.32180977 & 0.00410184 & 5.32E+02 \\
    5     & 363.4003476 & 341.5519319 & 2.94E+02 & 3.23E+02 & 364.9581233 & 0.36275628 & 0.00597385 & 3.23E+02 \\
    6     & 230.608119 & 207.585416 & 1.89E+02 & 1.89E+02 & 232.1016543 & 0.38903529 & 0.00799251 & 1.89E+02 \\
    7     & 147.8561494 & 124.2274595 & 1.18E+02 & 1.06E+02 & 149.3096521 & 0.40541155 & 0.0098586 & 1.06E+02 \\
    8     & 98.01126796 & 74.27302588 & 7.06E+01 & 5.63E+01 & 99.44065736 & 0.41527564 & 0.01127184 & 5.63E+01 \\
    9     & 69.13937735 & 45.768905 & 4.06E+01 & 2.79E+01 & 70.5547995 & 0.42098926 & 0.01205541 & 2.79E+01 \\
    10    & 53.13057813 & 30.57994962 & 2.21E+01 & 1.29E+01 & 54.53825576 & 0.42415733 & 0.01220761 & 1.29E+01 \\
    11    & 44.6597286 & 23.22992376 & 1.12E+01 & 5.69E+00 & 46.06330831 & 0.42583368 & 0.01190151 & 5.69E+00 \\
    12    & 40.41192388 & 20.07409231 & 5.29E+00 & 2.44E+00 & 41.81344865 & 0.4266743 & 0.01143556 & 2.44E+00 \\
    13    & 38.42463196 & 18.85806286 & 2.33E+00 & 1.02E+00 & 39.82519534 & 0.42706758 & 0.01106236 & 1.02E+00 \\
    14    & 37.56159752 & 18.41580876 & 9.70E-01 & 4.16E-01 & 38.9617434 & 0.42723837 & 0.01084908 & 4.16E-01 \\
    15    & 37.20752766 & 18.25657936 & 3.88E-01 & 1.65E-01 & 38.60750225 & 0.42730844 & 0.01074838 & 1.65E-01 \\
    16    & 37.06714943 & 18.19864095 & 1.52E-01 & 6.44E-02 & 38.46705611 & 0.42733622 & 0.01070538 & 6.44E-02 \\
    17    & 37.01251861 & 18.17727243 & 5.87E-02 & 2.49E-02 & 38.41239886 & 0.42734703 & 0.01068795 & 2.49E-02 \\
    18    & 36.99146859 & 18.16930635 & 2.25E-02 & 9.54E-03 & 38.39133866 & 0.42735119 & 0.01068108 & 9.54E-03 \\
    19    & 36.98340172 & 18.16631458 & 8.60E-03 & 3.65E-03 & 38.38326789 & 0.42735279 & 0.01067841 & 3.65E-03 \\
    20    & 36.98031974 & 18.16518558 & 3.28E-03 & 1.39E-03 & 38.38018441 & 0.4273534 & 0.01067738 & 1.39E-03 \\
    21    & 36.97914434 & 18.16475823 & 1.25E-03 & 5.31E-04 & 38.37900845 & 0.42735363 & 0.01067699 & 5.30E-04 \\
    22    & 36.97869653 & 18.16459616 & 4.76E-04 & 2.02E-04 & 38.37856042 & 0.42735372 & 0.01067684 & 2.02E-04 \\
    23    & 36.97852602 & 18.16453462 & 1.81E-04 & 7.69E-05 & 38.37838983 & 0.42735375 & 0.01067678 & 7.67E-05 \\
    24    & 36.97846112 & 18.16451124 & 6.90E-05 & 2.93E-05 & 38.3783249 & 0.42735377 & 0.01067676 & 2.93E-05 \\
    25    & 36.97843643 & 18.16450235 & 2.62E-05 & 1.11E-05 & 38.37830019 & 0.42735377 & 0.01067675 & 1.10E-05 \\
    26    & 36.97842703 & 18.16449898 & 9.99E-06 & 4.23E-06 & 38.37829079 & 0.42735377 & 0.01067675 & 4.12E-06 \\ \bottomrule
\caption{Iterations generated by the fractional quasi-Newton method.}
\end{longtable}
\end{footnotesize}

\item[ii)] $A_{g_f,\alpha_f} \Rightarrow \Phi \in \Ord_2^2(\xi)$ in $\Omega_\xi \subset B(\xi;\delta)$.

\begin{footnotesize}
\centering
\begin{longtable}{c|cccc|cccccccc}
\toprule
$i$& $[ x_i ]_2$ &$[ x_i ]_3$ &$\norm{x_i-x_{i-1}}_2$&$\norm{f_2(x_i)}_2$&$[ x_i ]_1$&$[ x_i ]_4$&$[ x_i ]_5$& $\norm{f_1(s_i)}_2$ \\
\midrule
    1     & 2029.854772 & 2022.247443 & 1.38E+03 & 2.00E+03 & 2032.218723 & 0.03297214 & 0.00056752 & 2.00E+03 \\
    2     & 1351.035349 & 1337.861649 & 9.64E+02 & 1.32E+03 & 1353.07091 & 0.16730757 & 0.00139631 & 1.32E+03 \\
    3     & 884.5286725 & 867.3584839 & 6.63E+02 & 8.49E+02 & 886.3385526 & 0.25962723 & 0.00256042 & 8.49E+02 \\
    4     & 570.3098992 & 550.3428213 & 4.46E+02 & 5.32E+02 & 571.9677708 & 0.32180977 & 0.00410184 & 5.32E+02 \\
    5     & 363.4003476 & 341.5519319 & 2.94E+02 & 3.23E+02 & 364.9581233 & 0.36275628 & 0.00597385 & 3.23E+02 \\
    6     & 230.608119 & 207.585416 & 1.89E+02 & 1.89E+02 & 232.1016543 & 0.38903529 & 0.00799251 & 1.89E+02 \\
    7     & 147.8561494 & 124.2274595 & 1.18E+02 & 1.06E+02 & 149.3096521 & 0.40541155 & 0.0098586 & 1.06E+02 \\
    8     & 98.01126796 & 74.27302588 & 7.06E+01 & 5.63E+01 & 99.44065736 & 0.41527564 & 0.01127184 & 5.63E+01 \\
    9     & 69.13937735 & 45.768905 & 4.06E+01 & 2.79E+01 & 70.5547995 & 0.42098926 & 0.01205541 & 2.79E+01 \\
    10    & 53.13057813 & 30.57994962 & 2.21E+01 & 1.29E+01 & 54.53825576 & 0.42415733 & 0.01220761 & 1.29E+01 \\
    11    & 36.97715447 & 18.16441312 & 2.04E+01 & 1.19E-03 & 38.37701761 & 0.42735403 & 0.0106761 & 1.19E-03 \\
    12    & 36.97842127 & 18.1644969 & 1.27E-03 & 7.75E-09 & 38.37828503 & 0.42735378 & 0.01067675 & 2.15E-07 \\
     \bottomrule
\caption{Iterations generated by the fractional quasi-Newton method accelerated.}
\end{longtable}
\end{footnotesize}

\end{itemize}
\end{example}

\begin{example}
Considering by hypothesis that $f_2\in \Ord_{NR}^2(\xi)$, and using the following values

\begin{eqnarray*}
DNI=94.3555,&
Tair=28.373,& \alpha=0.89964,
\end{eqnarray*}

the following iterations are obtained by using the fractional iterative methods given by the matrices  \eqref{eq:5-008} and \eqref{eq:5-009}.

\begin{itemize}
\item[i)] $A_{g_f,\beta} \Rightarrow \Phi \in \Ord_2^1(\xi)$.

\begin{footnotesize}
\centering
\begin{longtable}{c|cccc|cccccccc}
\toprule
$i$& $[ x_i ]_2$ &$[ x_i ]_3$ &$\norm{x_i-x_{i-1}}_2$&$\norm{f_2(x_i)}_2$&$[ x_i ]_1$&$[ x_i ]_4$&$[ x_i ]_5$& $\norm{f_1(s_i)}_2$ \\
\midrule
    1     & 2026.948258 & 2025.698601 & 1.38E+03 & 2.00E+03 & 2027.336246 & 0.03393789 & 0.00009324 & 2.00E+03 \\
    2     & 1346.157858 & 1343.993285 & 9.63E+02 & 1.32E+03 & 1346.49176 & 0.16860893 & 0.00022947 & 1.32E+03 \\
    3     & 878.348027 & 875.5257474 & 6.62E+02 & 8.47E+02 & 878.644763 & 0.26114907 & 0.00042095 & 8.47E+02 \\
    4     & 563.2981664 & 560.0145307 & 4.46E+02 & 5.31E+02 & 563.5698729 & 0.32347088 & 0.00067464 & 5.31E+02 \\
    5     & 355.8897756 & 352.2947261 & 2.94E+02 & 3.24E+02 & 356.1450042 & 0.36449952 & 0.00098289 & 3.24E+02 \\
    6     & 222.8349909 & 219.0449879 & 1.88E+02 & 1.90E+02 & 223.0796488 & 0.39081985 & 0.00131521 & 1.90E+02 \\
    7     & 139.9975258 & 136.1076542 & 1.17E+02 & 1.08E+02 & 140.2356025 & 0.4072064 & 0.00162172 & 1.08E+02 \\
    8     & 90.22025428 & 86.31555624 & 7.04E+01 & 5.77E+01 & 90.45437636 & 0.41705312 & 0.00185193 & 5.77E+01 \\
    9     & 61.57917913 & 57.74236258 & 4.05E+01 & 2.92E+01 & 61.81102578 & 0.42271878 & 0.00197603 & 2.92E+01 \\
    10    & 45.9860483 & 42.29258563 & 2.20E+01 & 1.37E+01 & 46.21665614 & 0.42580335 & 0.00199523 & 1.37E+01 \\
    11    & 38.07727079 & 34.56957756 & 1.11E+01 & 5.99E+00 & 38.3072503 & 0.42736783 & 0.00194279 & 5.99E+00 \\
    12    & 34.38320509 & 31.04574209 & 5.11E+00 & 2.46E+00 & 34.61289112 & 0.42809857 & 0.00187038 & 2.46E+00 \\
    13    & 32.785517 & 29.5629105 & 2.18E+00 & 9.76E-01 & 33.0150761 & 0.42841462 & 0.0018152 & 9.76E-01 \\
    14    & 32.13122049 & 28.97016965 & 8.83E-01 & 3.80E-01 & 32.36072761 & 0.42854405 & 0.00178421 & 3.80E-01 \\
    15    & 31.87136629 & 28.7388988 & 3.48E-01 & 1.47E-01 & 32.10085277 & 0.42859545 & 0.00176952 & 1.47E-01 \\
    16    & 31.7697037 & 28.64948568 & 1.35E-01 & 5.67E-02 & 31.9991821 & 0.42861556 & 0.00176316 & 5.67E-02 \\
    17    & 31.73020749 & 28.61501378 & 5.24E-02 & 2.19E-02 & 31.95968275 & 0.42862337 & 0.00176054 & 2.19E-02 \\
    18    & 31.71491342 & 28.60173066 & 2.03E-02 & 8.43E-03 & 31.94438747 & 0.4286264 & 0.00175949 & 8.43E-03 \\
    19    & 31.70900063 & 28.59661144 & 7.82E-03 & 3.25E-03 & 31.93847421 & 0.42862757 & 0.00175907 & 3.25E-03 \\
    20    & 31.70671661 & 28.59463794 & 3.02E-03 & 1.25E-03 & 31.93619001 & 0.42862802 & 0.00175891 & 1.25E-03 \\
    21    & 31.70583474 & 28.59387694 & 1.17E-03 & 4.84E-04 & 31.93530807 & 0.4286282 & 0.00175885 & 4.84E-04 \\
    22    & 31.70549433 & 28.59358344 & 4.50E-04 & 1.87E-04 & 31.93496763 & 0.42862826 & 0.00175882 & 1.87E-04 \\
    23    & 31.70536295 & 28.59347022 & 1.73E-04 & 7.20E-05 & 31.93483624 & 0.42862829 & 0.00175881 & 7.20E-05 \\
    24    & 31.70531225 & 28.59342654 & 6.69E-05 & 2.78E-05 & 31.93478553 & 0.4286283 & 0.00175881 & 2.78E-05 \\
    25    & 31.70529268 & 28.59340969 & 2.58E-05 & 1.07E-05 & 31.93476596 & 0.4286283 & 0.00175881 & 1.07E-05 \\
    26    & 31.70528513 & 28.59340319 & 9.96E-06 & 4.13E-06 & 31.93475841 & 0.4286283 & 0.00175881 & 4.13E-06 \\
    \bottomrule
\caption{Iterations generated by the fractional quasi-Newton method.}
\end{longtable}
\end{footnotesize}

\item[ii)] $A_{g_f,\alpha_f} \Rightarrow \Phi \in \Ord_2^2(\xi)$ in $\Omega_\xi \subset B(\xi;\delta)$.

\begin{footnotesize}
\centering
\begin{longtable}{c|cccc|cccccccc}
\toprule
$i$& $[ x_i ]_2$ &$[ x_i ]_3$ &$\norm{x_i-x_{i-1}}_2$&$\norm{f_2(x_i)}_2$&$[ x_i ]_1$&$[ x_i ]_4$&$[ x_i ]_5$& $\norm{f_1(s_i)}_2$ \\
\midrule
    1     & 2026.94826 & 2025.6986 & 1.38E+03 & 2.00E+03 & 2027.33625 & 0.03393789 & 0.0000932 & 2.00E+03 \\
    2     & 1346.15786 & 1343.99329 & 9.63E+02 & 1.32E+03 & 1346.49176 & 0.16860893 & 0.00022947 & 1.32E+03 \\
    3     & 878.348027 & 875.525747 & 6.62E+02 & 8.47E+02 & 878.644763 & 0.26114907 & 0.00042095 & 8.47E+02 \\
    4     & 563.298166 & 560.014531 & 4.46E+02 & 5.31E+02 & 563.569873 & 0.32347088 & 0.00067464 & 5.31E+02 \\
    5     & 355.889776 & 352.294726 & 2.94E+02 & 3.24E+02 & 356.145004 & 0.36449952 & 0.00098289 & 3.24E+02 \\
    6     & 222.834991 & 219.044988 & 1.88E+02 & 1.90E+02 & 223.079649 & 0.39081985 & 0.00131521 & 1.90E+02 \\
    7     & 139.997526 & 136.107654 & 1.17E+02 & 1.08E+02 & 140.235603 & 0.4072064 & 0.00162172 & 1.08E+02 \\
    8     & 90.2202543 & 86.3155562 & 7.04E+01 & 5.77E+01 & 90.4543764 & 0.41705312 & 0.00185193 & 5.77E+01 \\
    9     & 61.5791791 & 57.7423626 & 4.05E+01 & 2.92E+01 & 61.8110258 & 0.42271878 & 0.00197603 & 2.92E+01 \\
    10    & 45.9860483 & 42.2925856 & 2.20E+01 & 1.37E+01 & 46.2166561 & 0.42580335 & 0.00199523 & 1.37E+01 \\
    11    & 38.0772708 & 34.5695776 & 1.11E+01 & 5.99E+00 & 38.3072503 & 0.42736783 & 0.00194279 & 5.99E+00 \\
    12    & 31.7052694 & 28.5933984 & 8.74E+00 & 1.03E-05 & 31.9347427 & 0.42862831 & 0.0017588 & 1.03E-05 \\
    13    & 31.7052804 & 28.5933991 & 1.10E-05 & 2.80E-09 & 31.9347537 & 0.42862831 & 0.00175881 & 3.17E-08 \\
     \bottomrule
\caption{Iterations generated by the fractional quasi-Newton method accelerated.}
\end{longtable}
\end{footnotesize}

\end{itemize}
\end{example}

From the previous results it is observed that there is a considerable improvement in the order of convergence between the matrices \eqref{eq:5-008} and \eqref{eq:5-009}. Therefore, it may be established that it is more efficient to solve the function \eqref{eq:5-002} by implementing the fractional quasi-Newton method accelerated in the function \eqref{eq:5-005}. So, by providing multiple values of the parameters $DNI$ and $T_{air}$, it is possible to obtain a histogram of the efficiencies of a hybrid solar receiver analogous to the one shown in the Figure \ref{fig:02}. Finally, it is necessary to mention that the \mref{Corollary}{\ref{cor:4-001}} can also be implemented in the \textbf{generalized fractional quasi-Newton method}, which is obtained by using the matrix \eqref{eq:5-008} with the following function

\begin{eqnarray}
g_{a,b,f}(x):=af(x_i)+f^{(1)}(x_i)(x-bx_i), & a,b\in \nset{R},
\end{eqnarray}

as a consequence, denoting by $c$ an arbitrary constant, it is possible to define the following set of matrices

\begin{eqnarray}
\set{A_{g,\alpha} \ : \    A_{g,\alpha} \in  {}_n\IMa_{x,\alpha}^{1}(g) }\cap\set{A_{g,\alpha} \ : \  \left([A_{g,\alpha}]_{jk}(x) \right) :=\left( o_k^{\alpha}[g]_{j}(x)  \right)^{-1} }\cap \set{o^\alpha \ : \ o_k^{\alpha} c  \neq 0 \ \forall k\geq 1 } ,
\end{eqnarray}

and therefore, it is possible to define the following sets of fractional iterative methods

\begin{eqnarray}
\set{\Phi \ : \ A_{g,\beta}  \mbox{ uses }  g=g_{a,0,f} \mbox{ with }  0<a\leq 1}, \vspace{0.1cm}\\
\set{\Phi \ : \ A_{g,\beta} \mbox{ uses }  g=g_{1,b,f} \mbox{ with } 0< b \leq 1},
\end{eqnarray}

which correspond to two uncountable families of fractional fixed point methods in which the  \mref{Corollary}{\ref{cor:4-001}} can be implemented.

\section{Conclusions}

In all the examples shown, a decrease in the number of iterations necessary to converge to the solutions is observed when implementing the function \eqref{eq:4-005} in the fractional quasi-Newton method, which means that the generated sequences show an acceleration in their speed of convergence, which was to be expected given the \mref{Corollary}{\ref{cor:4-001}}. The fractional iterative methods, such as the fractional Newton-Raphson method, can find multiple zeros of a function using a single initial condition, this partially solves the intrinsic problem of classical iterative methods, which is that in general, to find $N$ zeros of a function, $N$ initial conditions must be provided. Due to the fractional operators implemented, these methods can be considered \textbf{non-local parametric iterative methods}, so they have two important characteristics:

\begin{itemize}
\item[i)] The initial condition does not necessarily have to be close to the sought values due to the non-local nature of fractional operators \cite{torres2020nonlinear}.

\item[ii)] When working in a space of $N$ dimensions, in the case that it is necessary to change the initial condition, unlike the classic iterative methods where in the worst case it is necessary to vary the $N$ components of the initial condition until obtaining a suitable value, in the fractional fixed point methods it is enough to vary the parameter $\alpha$ of the fractional operators until found an adequate value that allows generating a sequence that converges to a sought value \cite{torres2020fracsome}.
\end{itemize}

The above features make fractional iterative methods an ideal numerical tool for working with systems of nonlinear algebraic equations that vary with time-dependent parameters, as is the case of the functions \eqref{eq:5-002} and \eqref{eq:5-005}, which allows studying the behavior of temperatures and efficiencies of a hybrid solar receiver \cite{de2021fractional,torres2020reduction,torres2020fracpseunew}. 
Due to many nonlinear algebraic systems related to engineering and physics are often related to time-dependent parameters. Having a way of classifying and accelerating the order of convergence of fractional fixed point methods through the \mref{Corollary}{\ref{cor:4-001}} may become a fundamental piece to continue expanding the applications of the fractional operators.

\begin{figure}[!ht]
\centering
\includegraphics[width=0.6\textwidth, height=0.37\textwidth]{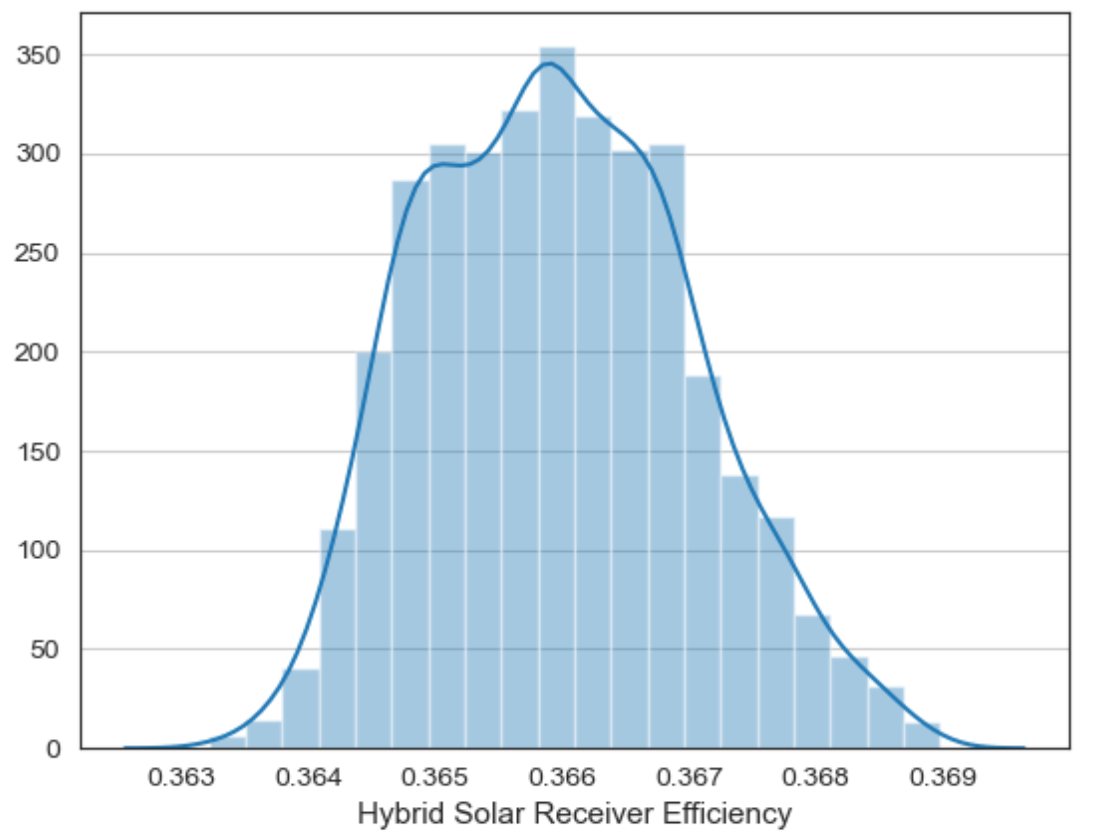}
\caption{Histogram and density curve of the efficiency of a hybrid solar receiver obtained from a simulation corresponding to a period of thirty days, which is equivalent to $2410$ pairs of parameters ($DNI,T_{air}$) randomly generated on the domain $[12,958]\times [11,45]$. The selected domain is based on data measured in real-time at the Center for Advanced Studies in Energy and Environment (CEAEMA) \cite{rodrigo2019performance,de2021fractional}. The values generated for the simulation presented the mean values $mean(DNI)=662.35$ and $mean(T_{air})=31.28$ with sample standard deviations $std(DNI)=257.83$ and $std(T_{air})=6.11
$, while the values of the efficiencies were obtained through the solutions of the function \eqref{eq:5-005} using the fractional quasi-Newton method accelerated.}\label{fig:02}
\end{figure}

\bibliography{Biblio}

\bibliographystyle{unsrt}

\end{document}